\documentclass[a4wide]{amsart}
\usepackage{amsmath}
\usepackage{amssymb}
\usepackage{latexsym}
\usepackage{amsbsy}
\usepackage[all]{xy}
\usepackage{amscd}
\usepackage[dvips]{graphicx}

\newtheorem{Prop}{Proposition}[section]
\newtheorem{Thm}[Prop]{Theorem}
\newtheorem{Lemma}[Prop]{Lemma}

\newtheorem{Remark}[Prop]{Remark}

\newtheorem{Definition}[Prop]{Definition}

    \def\bbq{{\mathbb Q}} \def\bb1{{\mathbb 1}}

\def\hom{\mbox{Hom}}
\def\aut{\mbox{Aut}}
\def\rad{\mbox{rad}\,}

\def\dim{\mbox{dim}\,}

\def\ed{\mbox{End}}

\def\mod{\mbox{mod}\,}  
    \def\aut{\mbox{Aut}\,}

\def\uq2{U_q(\hat{sl}_2)}

\def\bb{{\bf b}}

\def\nd{{\noindent}}

\def\mc{{\mathcal{C}}}

\def\md{{\mathcal{D}}}

\def\mh{{\mathcal{H}}}

\begin{document}
\title{Hall algebras associated to triangulated categories}

\thanks{ The research was
supported in part by NSF of China (No. 10631010) and by NKBRPC (No. 2006CB805905)  }

\author{Jie Xiao and Fan Xu }
\address{Department of Mathematical Sciences\\
Tsinghua University\\
Beijing 10084, P.~R.~China} \email{jxiao@math.tsinghua.edu.cn
(J.Xiao),\  f-xu04@mails.tsinghua.edu.cn (F.Xu)}

\maketitle

\begin{center}{\it \footnotesize Dedicated to  Yanan Lin on the occasion of his
50th birthday}
\end{center}
\bigskip

\bigskip

\begin{abstract}
 By counting with triangles and the octohedral axiom, we find a
direct way to prove the formula of To\"en in \cite{Toen2005}  for a
triangulated category with (left) homological-finite condition.
\end{abstract}
\section{Introduction}
Let $k$ be the finite field $\mathbb{F}_q$ with $q$ elements and
$A$ be a finite dimensional $k$-algebra.  Ringel associated the
module category of $A$ an associative algebra $\mh(A),$ which now
is called Ringel-Hall algebra, and used it to give a realization
of the positive part of simple Lie algebra when $A$ is  a
hereditary algebra of finite representation type (see
\cite{Ringel1990a} and \cite{Ringel1990b}). In general, the idea
of Ringel-Hall algebra constructs an assoicative algebra
$\mh(\mathcal{A})$ from an abelian category $\mathcal{A}.$ The
isomorphism classes of object in $\mathcal{A}$ generate the vector
space $\mathcal{H}(\mathcal{A})$ with multiplication
$[X]*[Y]=\sum_{[L]}F_{XY}^{L}[L],$ where $F_{XY}^{Z}$ is called
Hall number and is the number of subobject $L'$ of $L$ such that
$L'\cong X, L/L'\cong Y.$ Moreover, Ringel and Green showed when
$A$ is an arbitrary hereditary algebra, the composition subalgebra
of $\mh(A)$ gives a realization of the positive part of the
quantum enveloping algebra of the corresponding Kac-Moody
algebras( see \cite{Ringel1990} and \cite{green1995}, also
\cite{Ringel1996}). So the next question is, asked by Ringel in
\cite{Ringel1990}, to recover the whole Lie algebra and the whole
quantized enveloping algebra. A direct idea is to use Drinfeld
Double to piece together two Borel parts as showed in
\cite{Xiao1997}. However, this construction is not ``intrinsic",
i.e., not naturally induced by the module category of $A$.
Therefore, one need extend the module category of $A$ to a larger
category.

To deal with this question, several important developments have been
made. One is to use the 2-period triangulated category to define an
analog multiplication of Hall multiplication (see \cite{PX2000}).
Although this multiplication is not associative in general, the Lie
bracket induced by it satisfies Jacobi identity and a geometric
method  is verified feasible to define this Lie bracket directly
over the complex field (see \cite{XXZ2006}). On the other hand,
Kapranov in \cite{Kapranov1998} defined Heisenberg doubles for
hereditary categories and attached an associative algebra to the
derived category of any hereditary category. Recently, To\"en made a
remarkable development in this direction. He defined an associative
algebra corresponding to a dg category by using model categories and
fibre product of model categories. A key formula for the derived
Hall numbers, i.e. structure constant of the multiplication
analogous to Hall number, is given (see Proposition 5.1 in
\cite{Toen2005} or Section 3 in the following). In Remark 5.3 in
\cite{Toen2005}, To\"en asked whether one can define the derived
Hall algebra of any triangulated category under some finiteness
condition by his formula as the structure constant. The purpose of
this note is to find a direct method to deduce To\"en's formula for
arbitrary triangulated category under some finiteness conditions
(see section 2 for these conditions). We combine  the methods in
\cite{PX2000} and in \cite{Toen2005}. Our idea is as follows. First,
we recognize one of the main reasons that the multiplication defined
in \cite{PX2000} does not satisfy the associativity is that the
action of $\aut X\times \aut Y$ on $W(X,Y;L)$ is not free (see
Section 2 for the definitions of these and the following notations).
So we naturally consider the replacement $V(X,Y;L)$ by
$|\mathrm{Hom}(L,Y)_{X[1]}|/|\aut Y|$ or by
$|\mathrm{Hom}(X,L)_{Y}|/|\aut X|.$ Proposition \ref{4}$'$ shows the
explicit relation between $V(X,Y;L)$ and
$|\mathrm{Hom}(L,Y)_{X[1]}|/|\aut Y|$ or
$|\mathrm{Hom}(X,L)_{Y}|/|\aut X|.$ However, it is not proper to set
$|\mathrm{Hom}(L,Y)_{X[1]}|/|\aut Y|$ or
$|\mathrm{Hom}(X,L)_{Y}|/|\aut X|$ as the derived Hall numbers since
this definition is not symmetric, i.e.,
$|\mathrm{Hom}(L,Y)_{X[1]}|/|\aut Y|\neq
|\mathrm{Hom}(X,L)_{Y}|/|\aut X|.$ In fact, Proposition \ref{4}$'$
implies a symmetric expression.  A simple computation shows that the
proof of the associativity comes down to confirming the symmetry of
the other expression (see Definition \ref{def} and Proposition
\ref{5}), and the symmetry of the latter heavily depends on the
octahedral axiom.

Finally we should pay attention to the following points. One of the
next tasks is to construct To\"en's formula for a 2-period orbit
category. For the enveloping algebra $\mathcal{U}$ of a simple split
Lie algebra of type ADE. An arbitrarily large finite dimensional
quotient of $\mathcal{U}$ can be constructed in terms of
constructible functions on a ``triple variety" by Lusztig (see
\cite{Lusztig2000}) or in terms of the homology of a ``triple
variety" by Nakajima (see \cite{Nakajima1998}). It will be very
interesting to look for the relations between the construction by
``triple variety" and  the  To\"en's formula for a 2-period orbit
category.

\bigskip
\par\noindent {\bf Acknowledgments.}
The authors are very grateful to the referees for many helpful
comments. In particular, Proposition \ref{6} has been added and
some notations have been simplified for readability.

\section{Calculation with triangles}
Given a finite field $k$ with $q$ elements, let $\mc$ be a
$k$-additive triangulated category with the translation $T=[1].$ We
always assume in this paper that  $(1)$ the homomorphism space
$\hom_{\mc}(X,Y)$ for any two objects $X$ and $Y$ in $\mc$ is a
finite dimensional $k$-space, and $(2)$ the endomorphism ring $\ed
X$ for any indecomposable object $X$ is finite dimensional local
$k$-algebra. We note that the above two conditions imply the
Krull-Schmidt theorem holds in $\mc,$ i.e., any object in $\mc$ can
be decomposed into the direct sum of finitely many indecomposable
objects.  Moreover, we always assume that $\mc$ is (left) locally
homological finite, i.e., $\sum_{i\geq 0}\dim_k\hom(X[i],Y)<\infty$
for any $X$ and $Y$ in $\mc.$ We will use $fg$ to denote the
composition of morphisms $f: X\rightarrow Y$ and $g: Y\rightarrow
Z,$ and $|A|$ the cardinality of a finite set $A.$ For example, the
bounded derived category $\md^{b}(A)$ of the module category $\mod
A$ of a finite dimensional $k$-algebra $A$ satisfies all conditions
as above. However, its 2-period orbit category $\md^{b}(A)/T^{2},$
which is a triangulated category if $A$ is hereditary, does not
satisfy  the homological finiteness  property.

The following is an easy result  in \cite{PX2000}.
\begin{Lemma}
Given a triangle of form
\begin{equation}
\xymatrix{M\ar[rr]^-{(\begin{array}{cc}
  f_1 & f_2 \\
\end{array})}&& N_1\oplus N_2\ar[rr]^-{\left(%
\begin{array}{c}
  g_1 \\
  g_2 \\
\end{array}%
\right)}&& L\ar[rr]^{\begin{array}{c}
  h \\
\end{array}}&& M[1]}
\end{equation} If $f_2=0,$ then it
is isomorphic to the triangle of form
\begin{equation}
\xymatrix{M\ar[rr]^-{( \begin{array}{cc}
  f_1 & f_2 \\
\end{array})}&& N_1\oplus N_2\ar[rr]^-{\left(%
\begin{array}{cc}
  g_{11} & 0 \\
  0 & g_{22} \\
\end{array}%
\right)}&& L_1\oplus L_2\ar[rr]^{\left(%
\begin{array}{c}
  h_1 \\
  0 \\
\end{array}%
\right)}&& M[1]}
\end{equation}
where $g_{22}: N_2\rightarrow L_2$ is an isomorphism and
\begin{equation}
\xymatrix{M\ar[r]^{f_1}& N_1\ar[r]^{g_{11}}& L_1\ar[r]^{h_1}&
M[1]} \end{equation} is a triangle.
\end{Lemma}
The lemma shows the triangle $(1)$ is the direct sum of $(3)$ and
the following contractible triangle (also see \cite{Neeman}):
$$
\xymatrix{0\ar[r]& N_2\ar[r]^{\sim}& L_2\ar[r]& 0}
$$
Now we recall some notation in \cite{PX2000}.

Given any object $X,Y,Z,L,L'$ and $M$ in $\mc,$ we define
$$W(X,Y;L)=\{(f,g,h)\in \hom(X,L)\times
\hom(L,Y)\times(Y,X[1])\mid$$$$
X\xrightarrow{f}L\xrightarrow{g}Y\xrightarrow{h}X[1] \mbox{ is a
triangle}\}$$ The following we simply write $(f,g,h)$ is a
triangle to denote
$X\xrightarrow{f}L\xrightarrow{g}Y\xrightarrow{h}X[1] \mbox{ is a
triangle}.$

The action of $\aut X\times \aut Y$ on $W(X,Y;L)$ induces the
orbit space
$$V(X,Y;L)=\{(f,g,h)^{\wedge}\mid (f,g,h)\in W(X,Y;L)\}$$ where
$$(f,g,h)^{\wedge}=\{(af,gc^{-1},ch(a[1])^{-1})\mid (a,c)\in \aut X\times \aut Y\}$$
For any $((f,g,h)(l,m,n))\in W(X,Y;L)\times W(Z,L;M),$ we define
$$
\hspace{-9cm}((f,g,h),(l,m,n))^{*}$$$$=\{((fb^{-1},bgc^{-1},ch),(dl,mb^{-1},bn(d[1])^{-1}))\mid
b\in \aut L,c\in \aut Y,d\in \aut Z\}
$$
So we have the orbit space
$$
\hspace{-7cm}(W(X,Y;L)\times W(Z,L;M))^{*}$$
$$=\{((f,g,h),(l,m,n))^{*}\mid
((f,g,h),(l,m,n))\in W(X,Y;L)\times W(Z,L;M)\}
$$
Dually, for any $((l',m',n'),(f',g',h'))\in W(Z,X;L')\times
W(L',Y;M)$ we define
$$
\hspace{-8cm}((l',m',n'),(f',g',h'))^{*}$$$$=\{((dl'b'^{-1},b'm',n'(d[1])^{-1}),(b'f',g'c^{-1},ch'(b'[1])^{-1}))\mid
d\in \aut Z,b'\in \aut L',c\in \aut Y\}
$$
Then the orbit space is
$$ \hspace{-7cm}(W(Z,X;L')\times
W(L',Y;M))^{*}$$
$$=\{((l',m',n'),(f',g',h'))^{*}\mid
((l',m',n'),(f',g',h'))\in W(X,Y;L)\times W(Z,L;M)\}
$$

We define the action of $\aut Z$ on $W(Z,L;M)$ as follows: For any
$(l,m,n)\in W(Z,L;M),$ $d(l,m,n)=(dl,m,n(d[1])^{-1})$ for any
$(l,m,n)\in W(Z,L;M)$ and any $d\in \aut Z.$ We denote the orbit
by
$$
(l,m,n)^{*}_{Z}=\{(dl,m,n(d[1])^{-1})\mid d\in \aut Z\}
$$
then the orbit space is
$$
W(Z,L;M)^{*}_{Z}=\{(l,m,n)^{*}_{Z}\mid (l,m,n)\in W(Z,L;M)\}.
$$
Dually, we also have the action of $\aut L$ on $W(Z,L;M).$
$$(l,m,n)^{*}_{L}=\{(l,mb^{-1},bn)\mid b\in \aut L\}$$
and
$$
W(Z,L;M)^{*}_{L}=\{(l,m,n)^{*}_{L}\mid (l,m,n)\in W(Z,L;M)\}.$$ We
have the following diagram to help understanding (see
\cite{PX2000}).

\begin{equation}
\xymatrix{Z\ar@{=}[r]\ar[d]^{l'} &Z\ar[d]^l\\
L'\ar@{.>}[r]^{f'}\ar[d]^{m'} &M\ar@{.>}[r]^{g'}\ar[d]^m &Y\ar@{.>}[r]^-{h'}\ar@{=}[d] &L'[1]\ar[d]^{m'[1]}\\
X\ar[r]^f\ar[d]^{n'} &L\ar[r]^g\ar[d]^n &Y\ar[r]^-{h} &X[1]\\
Z[1]\ar@{=}[r] &Z[1]}
\end{equation}

Let $X,Y\in\mc.$ Define $\rad\hom(X,Y),$ the radical of
$\hom_{\mc}(X,Y),$ to be

$$\nd \rad\hom(X,Y)=\{f\in \hom_{\mc}(X,Y)\mid gfh \mbox{ is not an
isomorphism }$$ $$\mbox{for any } g:A\rightarrow X\mbox{ and
}h:Y\rightarrow A \mbox{ with indecomposable } A \}$$
\begin{Lemma}
Let $X,Y\in\mc$ and $n\in \mathrm{Hom}_{\mc}(X,Y),$  then there
exists the decompositions $X=X_1(n)\oplus X_2(n)$,\
$Y=Y_1(n)\oplus Y_2(n)$ and $a\in \aut X,$ $c\in \aut Y$ such that
$a n c=\left(
         \begin{array}{cc}
           n'_{11} & 0 \\
           0 & n'_{22} \\
         \end{array}
       \right),$  where $n'_{11}$ is an isomorphism between $X_1(n)$ and $Y_1(n),$   $n'_{22}\in
       \mathrm{radHom}(X_2(n),Y_2(n)).$
\end{Lemma}
\begin{proof}
Let $X=\bigoplus_{i} X_i$ and $Y=\bigoplus_{j}Y_j$ be the direct
sums of  indecomposable objects. For any indecomposable summands
$X_i$ of $X$ and $Y_j$ of $Y,$ the morphism $n$ induces the
morphism $n_{ij}\in \hom(X_i,Y_j).$ If $n_{ij}\in
\rad\hom(X_i,Y_j)$ for all $i,j,$ then $n\in \rad\hom(X,Y).$ So we
only need to take $X_1=Y_1=0.$ If there exist some isomorphism
$n_{ij},$ we may assume it is $n_{11}$ without loss of generality,
then we have
$$
\xymatrix{X=X_1\oplus X'\ar[rr]^{\left(
         \begin{array}{cc}
           n_{11} & n_{12} \\
           n_{21} & n_{22} \\
         \end{array}
       \right)}& & Y=Y_1\oplus Y'}
$$
Consider $\left(
         \begin{array}{cc}
           1 & 0 \\
           -n_{21}n^{-1}_{11} & 1 \\
         \end{array}
       \right)\in \aut X$ and $\left(
         \begin{array}{cc}
           1 & -n^{-1}_{11}n_{12} \\
           0 & 1 \\
         \end{array}
       \right)\in \aut Y,$ then
$$
\left(
         \begin{array}{cc}
           1 & 0 \\
           -n_{21}n^{-1}_{11} & 1 \\
         \end{array}
       \right) \left(
         \begin{array}{cc}
           n_{11} & n_{12} \\
           n_{21} & n_{22} \\
         \end{array}
       \right)\left(
         \begin{array}{cc}
           1 & -n^{-1}_{11}n_{12} \\
           0 & 1 \\
         \end{array}
       \right)=\left(
         \begin{array}{cc}
           n_{11} & 0 \\
           0 & n'_{22} \\
         \end{array}
       \right).
$$
The similar discussion works for $n'_{22}.$  By induction, we
achieve the claim of the lemma.
\end{proof}
\begin{Remark}
Any $(l,m,n)^{\wedge}\in V(Z,L;M)$ has the representative of the
form:
\begin{equation}
\xymatrix{Z\ar[rr]^{\left(%
\begin{array}{c}
  0 \\
  l_2 \\
\end{array}%
\right)}&& M\ar[rr]^{(%
\begin{array}{cc}
  0 & m_2 \\
\end{array}%
)}&& L\ar[rr]^{\left(%
\begin{array}{cc}
  n_{11} & 0 \\
  0 & n_{22} \\
\end{array}%
\right)}&& Z[1]}
\end{equation}
where $Z=Z_1\oplus Z_2,$ $L=L_1\oplus L_2$ and $n_{11}$ is an
isomorphism between $L_1$ and $Z_1[1],$ $n_{22}\in
\mathrm{radHom}(L_2,Z_2[1]).$ Depending on Lemma 1.1, it can be
decomposed into the direct sum of two triangles, in which one is a
contractible triangle.
\end{Remark}
In order to simplify the notation, for $X,Y\in\mc,$ we set
$$\{X,Y\}=\prod_{i>0}|\hom(X[i],Y)|^{(-1)^{i}}.$$
\begin{Lemma}\label{1}
For any $\alpha=(l,m,n)\in W(Z,L;M),$ set
$$
n\mathrm{Hom}(Z[1],L)=\{b\in \mathrm{End}L\mid  b=nt\ \mbox{for
some}\ t\in \mathrm{Hom}(Z[1],L)\}
$$
and
$$\mathrm{Hom}(Z[1],L)n=\{d\in \mathrm{End}Z[1]\mid d=sn \ \mbox{for some}\ s\in \mathrm{Hom}(Z[1],L) \} $$
We have

\begin{itemize}
\item[(1)]
 $|n\mathrm{Hom}(Z[1],L)| =
\prod_{i>0}\frac{|\mathrm{Hom}(M[i],L)|^{(-1)^{i}}}{|\mathrm{Hom}(Z[i],L)|^{(-1)^{i}}|\mathrm{Hom}(L[i],L)|^{(-1)^{i}}}
=\frac{\{M,L\}}{\{Z,L\}\{L,L\}}$
\item[(2)]$|\mathrm{Hom}(Z[1],L)n| =
\prod_{i>0}\frac{|\mathrm{Hom}(Z[i],M)|^{(-1)^{i}}}{|\mathrm{Hom}(Z[i],L)|^{(-1)^{i}}|\mathrm{Hom}(Z[i],Z)|^{(-1)^{i}}}=\frac{\{Z,M\}}{\{Z,L\}\{Z,Z\}}$

\end{itemize}

\end{Lemma}
\begin{proof}We only prove the first identity. It is similar to
prove the second. Applying $\hom(-,L)$ to the triangle
$Z\xrightarrow{l}M\xrightarrow{m}L\xrightarrow{n}Z[1]$ we get the
long exact sequence
$$
\cdots\rightarrow \hom(L[1],L)\rightarrow \hom(M[1],L)\rightarrow
\hom(Z[1],L)\xrightarrow{n^{*}}\hom(L,L)\rightarrow \cdots
$$
Since $n\hom(Z[1],L)=\mbox{Image of }n^{*},$ we have the identity
in this lemma.
\end{proof}

By Lemma 1.2, for any $L\xrightarrow{n}Z[1],$  there exist the
decompositions $L=L_{1}(n)\oplus L_{2}(n)$,
$Z[1]=Z_{1}[1](n)\oplus Z_2[1](n)$ and $b\in \aut L,$ $d\in \aut
Z$ such that $bn(d[1])^{-1}=\left(
         \begin{array}{cc}
           n_{11} & 0 \\
           0 & n_{22} \\
         \end{array}
       \right)$ and the induced maps
$n_{11}:L_{1}(n)\rightarrow Z_{1}[1](n)$ is an isomorphism and
$n_{22}:L_{2}(n)\rightarrow Z_2[1](n)$ contains no isomorphism
component. The above decomposition only depends on the equivalence
class  of $n$ up to an isomorphism. Let $\alpha=(l,m,n)^{\wedge}\in
V(Z,L;M),$ the classes of $\alpha$ and $n$ are determined to each
other in $V(Z,L;M).$ We may  denote $n$ by $n(\alpha)$ and $L_1(n)$
by $L_1(\alpha)$ respectively.

The following is a refinement of Lemma 7.1 in \cite{PX2000}.
\begin{Prop}\label{4}  We have
$$
|W(Z,L;M)^{*}_{Z}|=\sum_{\alpha\in V(Z,L;M)}\frac{|\mathrm{Aut}
L||\mathrm{End}
L_1(\alpha)|}{|n(\alpha)\mathrm{Hom}(Z[1],L)||\mathrm{Aut}
L_1(\alpha)|}
$$
\end{Prop}
\begin{proof}
Consider the map
$$
\varphi: W(Z,L;M)^{*}_{Z}\rightarrow V(Z,L;M)
$$
sending $(l,m,n)^{*}$ to $(l,m,n)^{\wedge}.$ The action of $\aut
L$ on $W(Z,L;M)$ naturally induces the action on
$W(Z,L;M)^{*}_{Z}$ with the orbit space $V(Z,L;M).$ Hence, for any
$(l,m,n)^{\wedge}\in V(Z,L;M),$ $\varphi^{-1}((l,m,n)^{\wedge})$
is just the orbit of $(l,m,n)^{*}$ under the action of $\aut L.$
Let $\alpha=(l,m,n)^{\wedge}.$ For $(l,m,n)^{*},$ we denote its
stable subgroup by $G((l,m,n)^{*}).$ Let $ n'=bn(d[1])^{-1}=\left(
         \begin{array}{cc}
           n_{11} & 0 \\
           0 & n_{22} \\
         \end{array}
       \right), l'=dl$ and $m'=mb^{-1}$ where $b\in \aut L$ and $d\in \aut Z.$  Then
       $$|G((l,m,n)^{*})|=|G((l',m',n')^{*})|
$$
$$
G((l',m',n')^{*})=\{b\in \aut L\mid
(l',m'b^{-1},bn')^{*}=(l',m',n')^{*}\}
$$
By definition, we have $G((l',m',n')^{*})=\{b\in \aut L\mid
m'b=m'\},$ i.e.
$$1-G((l',m',n')^{*})=\{b'\in \ed L\mid m'b'=0 \mbox{ and
}1-b'\in \aut L\}$$ On the other hand, $m'b'=0$ if and only if
$b'=n't$ for some $t\in \hom(Z[1],L).$ Let $t=\left(
         \begin{array}{cc}
           t_{11} & t_{12} \\
           t_{21} & t_{22} \\
         \end{array}
       \right)$ with respect to the above decompositions. Then we
       need
$$b=\left(
         \begin{array}{cc}
           1-n_{11}t_{11} & -n_{11}t_{12} \\
           -n_{22}t_{21} & 1-n_{22}t_{22} \\
         \end{array}
       \right)\in \aut L$$
The morphism $n_{22}\in \rad\hom(L_2,Z_2[1])$ implies
$n_{22}t_{21}$ is nilpotent and $1-n_{22}t_{22}\in \aut L_{2}.$
Hence, $b\in \aut L$ if and only if $1-n_{11}t_{11}\in \aut L_1.$
We obtain
$$
\hspace{-10cm}|G((l',m',n')^{*})|
$$
$$=|\aut L_1(\alpha)||n_{11}\hom(Z_1[1](\alpha),L_2(\alpha))||n_{22}\hom(Z_2[1](\alpha),L_1(\alpha))||n_{22}\hom(Z_2[1],L_2)|
$$
$$
\hspace{-6.2cm}=|n\hom(Z[1],L)||\aut L_1(\alpha)|/|\ed
L_1(\alpha)|
$$
On the other hand, we have $|n\hom(Z[1],L)|=|n'\hom(Z[1],L)|$ by
the above remark. We complete the proof of the proposition.
\end{proof}
Dually, we also have
$$
|W(Z,L;M)^{*}_{L}|=\sum_{\alpha\in V(Z,L;M)}\frac{|\aut Z||\ed
Z_1(\alpha)|}{|\hom(Z,L[-1])n(\alpha)[-1]||\aut Z_1(\alpha)|}
$$
where $Z_1(\alpha)\simeq L_{1}(\alpha)[-1].$

 As in \cite{Toen2005}, we denote by $\mathrm{Hom}(Z,M)_{L}$ the subset of $\hom(Z,M)$ consisting of
morphisms $l:Z\rightarrow M$ whose cone $Cone(l)$ is isomorphic to
$L.$ There is a natural action of the group $\aut Z$ on
$\mathrm{Hom}(Z,M)_{L}$ by $d\cdot l=dl,$ the orbit is denoted by
$l^{*}$ and the orbit space is denoted by
$\mathrm{Hom}(Z,M)^{*}_{L}.$ Dually We also have the subset
$\mathrm{Hom}(M,L)_{Z[1]}$ of $\hom(M,L)$ with the group action of
$\aut L$ and the orbit space $\mathrm{Hom}(M,L)^{*}_{Z[1]}.$
\begin{Prop}\label{abc}
There exist bijections:
$$
W(Z,L;M)^{*}_{Z}\rightarrow \mathrm{Hom}(M,L)_{Z[1]}\mbox{ and
}W(Z,L;M)^{*}_{L}\rightarrow \mathrm{Hom}(Z,M)_{L}.
$$
Moreover,
$|V(Z,L;M)|=|\mathrm{Hom}(M,L)^{*}_{Z[1]}|=|\mathrm{Hom}(Z,M)^{*}_{L}|.$
\end{Prop}
\begin{proof}
We have the natural surjections:
$$
\xymatrix{W(Z,L;M)^{*}_{Z}&W(Z,L;M)\ar[l]_-{
\pi_1}\ar[r]^-{\pi_2}& \mathrm{Hom}(M,L)_{Z[1]}}
$$
$\pi^{-1}_1((l,m,n)^{*})=\{(dl,m,n(d[1])^{-1})\mid d\in \aut Z\}$
and\\
$\pi^{-1}_2(m)=\{(l,m,n)\mid (l,m,n) \mbox{ is a triangle}\}.$ It
is clear that $\pi^{-1}_1((l,m,n)^{*})\subseteq \pi^{-1}_2(m).$ On
the other hand, for any $(l_1,m,n_1),(l_2,m,n_2)\in
\pi^{-1}_2(m),$ We have the following diagram:
$$
\xymatrix{Z\ar[r]^{l_1}& M\ar[r]^{m}\ar@{=}[d]&
L\ar[r]^{n_1}\ar@{=}[d]&Z[1]\\
Z\ar[r]^{l_2}\ar@{-->}[u]^{d}& M\ar[r]^{m}& L\ar[r]^{n_2}&Z[1]}
$$
Using (Tr3) in the triangulated category axioms, there exists an
isomorphism $d\in \mathrm{Aut}Z$ such that the above diagram
commutative. So $(l_2,m,n_2)=(dl_1,m,n_1(d[1])^{-1}).$ This shows
$\pi^{-1}_1((l,m,n)^{*})= \pi^{-1}_2(m).$ Hence, we naturally
define a bijection: $W(Z,L;M)^{*}_{Z}\rightarrow
\mathrm{Hom}(M,L)_{Z[1]}$ sending $(l,m,n)^{*}$ to $m.$ Dually, we
also have a bijection: $W(Z,L;M)^{*}_{L}\rightarrow
\mathrm{Hom}(Z,M)_{L}.$ Consider the maps:
$$
\varphi: W(Z,L;M)^{*}_{Z}\rightarrow V(Z,L;M)
$$
sending $(l,m,n)^{*}$  to $(l,m,n)^{\wedge}$ and
$$
\phi: \mathrm{Hom}(M,L)_{Z[1]}\rightarrow
\mathrm{Hom}(M,L)_{Z[1]}^{*}
$$
sending $m$ to $m^{*}.$ Then
$$
\varphi^{-1}((l,m,n)^{\wedge})=\{(l,mb^{-1},bn)^{*}\mid b\in \aut
L\}
$$
and
$$
\phi^{-1}(m^{*})=\{mb^{-1}\mid b\in \aut L\}
$$
the bijection between $W(Z,L;M)^{*}_{Z}$ and
$\mathrm{Hom}(M,L)_{Z[1]}$ induces the bijection between
$\varphi^{-1}((l,m,n)^{\wedge})$ and $\phi^{-1}(m^{*}).$ This
shows $|V(Z,L;M)|=|\mathrm{Hom}(M,L)^{*}_{Z[1]}|.$
\end{proof}
By Proposition \ref{abc}, Proposition \ref{4} and its dual formula
can be rewritten as follows.

\bigskip
\nd {\bf Proposition 2.5$'$} The following equalities hold.
$$
\frac{|\mathrm{Hom}(M,L)_{Z[1]}|}{|\mathrm{Aut}L|}\cdot
\frac{\{M,L\}}{\{Z,L\}\cdot\{L,L\}}=\sum_{\alpha\in
V(Z,L;M)}\frac{|\mathrm{End}L_1(\alpha)|}{|\mathrm{Aut}L_1(\alpha)|}
$$
$$
\frac{|\mathrm{Hom}(Z,M)_{L}|}{|\mathrm{Aut}Z|}\cdot
\frac{\{Z,M\}}{\{Z,L\}\cdot\{Z,Z\}}=\sum_{\alpha\in
V(Z,L;M)}\frac{|\mathrm{End}L_1(\alpha)|}{|\mathrm{Aut}L_1(\alpha)|}
$$

\section{Hall algebra arising in a triangulated category}
Let $\mc$ be a $k$-additive and $\mc$ a (left) locally homological
finite triangulated category. For any $X,Y,L\in \mc,$ by
 Proposition 2.5$'$, we define
\begin{eqnarray}
  F_{XY}^{L}:&=& \{X,Y\}\sum_{\alpha\in
V(X,Y;L)}\frac{|\ed X_1(\alpha)|}{|\aut X_1(\alpha)|} \nonumber \\
   &=&
   \frac{|\mathrm{Hom}(L,Y)_{X[1]}|}{|\mathrm{Aut}(Y)|\{Y,Y\}}\{L,Y\}=\frac{|\mathrm{Hom}(X,L)_{Y}|}{|\mathrm{Aut}(X)|\{X,X\}}\{X,L\}\nonumber
\end{eqnarray}
This formula is called To\"en's formula (\cite[Proposition
5.1]{Toen2005}). We will define an associative algebra arising from
$\mc,$ by using $F_{XY}^{L}$ as structure constants. For any $X\in
\mc,$ we denote its isomorphism class by $[X].$ Let $\mh$ be the
$\bbq$-space with the basis $\{u_{[X]}\mid X\in\mc\}.$ We define
$$
u_{[X]}*u_{[Y]}=\sum_{[L]}F_{XY}^{L}u_{[L]}
$$
Since $\hom(Y,X[1])$ is a finite set, the sum only has finitely many
nonzero summands.
\begin{Definition}\label{def}
Given any objects $X,M$ and $L,L'$ in $\mc,$ we set
$$
W(L',L;M\oplus X)=\{((f',-m'),\left(%
\begin{array}{c}
  m \\
  f \\
\end{array}%
\right),\theta)\mid ((f',-m'),\left(%
\begin{array}{c}
  m \\
  f \\
\end{array}%
\right),\theta)\mbox{ is a triangle }\}.
$$
For fixed $Y,Z\in\mc,$ we define its subsets
$$
W(L',L;M\oplus X)^{Y,Z}_{L}=\{((f',-m'),\left(%
\begin{array}{c}
  m \\
  f \\
\end{array}%
\right),\theta)\mid ((f',-m'),\left(%
\begin{array}{c}
  m \\
  f \\
\end{array}%
\right),\theta)\mbox{ is a triangle, }$$$$Cone(m)\cong Z[1],
\mbox{and }Cone(f)\cong Y\}.
$$
and
$$
W(L',L;M\oplus X)^{Y,Z}_{L'}=\{((f',-m'),\left(%
\begin{array}{c}
  m \\
  f \\
\end{array}%
\right),\theta)\mid ((f',-m'),\left(%
\begin{array}{c}
  m \\
  f \\
\end{array}%
\right),\theta)\mbox{ is a triangle, }$$$$Cone(m')\cong Z[1],
\mbox{and }Cone(f')\cong Y\}.
$$
\end{Definition}

In fact, we have
\begin{Prop}\label{6} The equality holds.
$$
W(L',L;M\oplus X)^{Y,Z}_{L}=W(L',L;M\oplus X)^{Y,Z}_{L'}
$$
\end{Prop}
The proof of Proposition~\ref{6} needs  the octahedral axiom and
pushout property in triangulated category. The following property
can be founded in \cite{PS1988}.
\begin{Prop}\label{7}
The following condition is equivalent to the octahedral axiom:
\\
Given a ``pushout'' square,i.e., a commutative square
\begin{equation}
\xymatrix{L'\ar[d]^{m'}\ar[r]^{f'}&M\ar[d]^{m}\\X\ar[r]^{f}&L}
\end{equation}
forming a distinguished triangle
$$
\xymatrix{L'\ar[rr]^-{(%
\begin{array}{cc}
  f' & -m' \\
\end{array}%
)}&&M\oplus X\ar[rr]^-{\left(%
\begin{array}{c}
  m \\
  f \\
\end{array}%
\right)}&&L\ar[rr]^{\begin{array}{c}
  \theta \\
\end{array}}&&L'[1]}
$$
and $\theta$ as above, it can be extended to a commutative diagram
\begin{equation}
\xymatrix{L'\ar[r]^{f'}\ar[d]^{m'} &M\ar@{.>}[r]^{g'}\ar[d]^m &Y\ar@{.>}[r]^-{h'}\ar@{=}[d] &L'[1]\ar[d]^{m'[1]}\\
X\ar[r]^f &L\ar[r]^g&Y\ar[r]^-{h} &X[1]}
\end{equation}
with $\theta=-gh'$.
\end{Prop}

\vspace{.2in}

\noindent \textit{Proof of Proposition~\ref{6}.} For any $((f',-m'),\left(%
\begin{array}{c}
  m \\
  f \\
\end{array}%
\right),\theta)\in W(L',L;M\oplus X)^{Y,Z}_{L},$ we have the
 corresponding diagram as follows:
\begin{equation}\label{8}
\xymatrix{&Z\ar[d]^l\\
L'\ar[r]^{f'}\ar[d]^{m'} &M \ar[d]^m &&\\
X\ar[r]^f &L\ar[r]^g\ar[d]^n &Y\ar[r]^-{h} &X[1]\\
 &Z[1]}
\end{equation}

Using Proposition~\ref{7} in horizontal and vertical direction two
times, diagram~(\ref{8}) is changed into:
\begin{equation}
\xymatrix{Z\ar@{=}[r]\ar@{.>}[d]^{l'} &Z\ar[d]^l\\
L'\ar[r]^{f'}\ar[d]^{m'} &M\ar@{.>}[r]^{g'}\ar[d]^m &Y\ar@{.>}[r]^-{h'}\ar@{=}[d] &L'[1]\ar[d]^{m'[1]}\\
X\ar[r]^f\ar@{.>}[d]^{n'} &L\ar[r]^g\ar[d]^n &Y\ar[r]^-{h} &X[1]\\
Z[1]\ar@{=}[r] &Z[1]}
\end{equation}
This shows $((f',-m'),\left(%
\begin{array}{c}
  m \\
  f \\
\end{array}%
\right),\theta)\in W(L',L;M\oplus X)^{Y,Z}_{L'},$ so
$W(L',L;M\oplus X)^{Y,Z}_{L}\subseteq W(L',L;M\oplus
X)^{Y,Z}_{L'}.$ Similarly, $W(L',L;M\oplus X)^{Y,Z}_{L'}\subseteq
W(L',L;M\oplus X)^{Y,Z}_{L}.$\qed

  Now we can set
  $$W_{Y,Z}(L',L;M\oplus X):=W(L',L;M\oplus X)^{Y,Z}_{L'}=
W(L',L;M\oplus X)^{Y,Z}_{L}$$
Define
$$
\mathrm{Hom}(M\oplus X,L)^{Y,Z[1]}_{L'[1]}=\{(m,f)\in
\mathrm{Hom}(M\oplus X,L)\mid$$$$\hspace{3cm} Cone(f)\simeq Y,
Cone(m)\simeq Z[1] \mbox{ and }Cone(m,f)\simeq L'[1]\}
$$
and
$$
\mathrm{Hom}(L',M\oplus X)^{Y,Z[1]}_{L}=\{(f',-m')\in
\mathrm{Hom}(L',M\oplus X)\mid $$$$\hspace{3cm}Cone(f')\simeq Y,
Cone(m')\simeq Z[1] \mbox{ and } Cone(f',-m')\simeq L\}
$$

The group  actions of $\aut L'$ and $\aut L$ on $W(L',L;M\oplus
X)$ naturally induce the actions on $W_{Y,Z}(L',L;M\oplus X).$ By
Proposition \ref{abc}, the orbit spaces under the action of $\aut
L'$ and $\aut L$ are $\mathrm{Hom}(M\oplus X,L)^{Y,Z[1]}_{L'[1]}$
and $\mathrm{Hom}(L',M\oplus X)^{Y,Z[1]}_{L},$ respectively. Under
the group action of $\aut L\times \aut L',$ the orbit space of
$W_{Y,Z}(L',L;M\oplus X)$ is denoted by $V_{Y,Z}(L',L;M\oplus X).$
Of course, $V_{Y,Z}(L',L;M\oplus X)$ is a subset of
$V(L',L;M\oplus X).$ Naturally, we have the following commutative
diagram:
$$
\xymatrix{W_{Y,Z}(L',L;M\oplus X)\ar[r]\ar[d]&\mathrm{Hom}(M\oplus X,L)^{Y,Z[1]}_{L'[1]}\ar[d]\\
\mathrm{Hom}(L',M\oplus X)^{Y,Z[1]}_{L}\ar[r]&V_{Y,Z}(L',L;M\oplus
X)}
$$
Applying Proposition \ref{4}$'$, we have
\begin{Prop}\label{5} The equalities hold.
$$
\frac{|\mathrm{Hom}(M\oplus X,L)^{Y,Z[1]}_{L'[1]}|}{|\mathrm{Aut}
L|}\frac{\{M\oplus X,L\}}{\{L',L\}\{L,L\}}=\sum_{\alpha\in
V_{Y,Z}(L',L;M\oplus X)}\frac{|\mathrm{End}
L_1(\alpha)|}{|\mathrm{Aut} L_1(\alpha)|}.
$$
and
$$
\frac{|\mathrm{Hom}(L',M\oplus X)^{Y,Z[1]}_{L}|}{|\mathrm{Aut}
L'|}\frac{\{L',M\oplus X\}}{\{L',L\}\{L',L'\}}=\sum_{\alpha\in
V_{Y,Z}(L',L;M\oplus X)}\frac{|\mathrm{End}
L_1(\alpha)|}{|\mathrm{Aut} L_1(\alpha)|}.
$$
\end{Prop}

\begin{Prop}\label{9}
There exist bijections:
$$
\mathrm{Hom}(X,L)_{Y}\times \mathrm{Hom}(M,L)_{Z[1]}\rightarrow
\bigcup_{[L']}\mathrm{Hom}(M\oplus X,L)^{Y,Z[1]}_{L'[1]}
$$
and
$$
\mathrm{Hom}(L',X)_{Z[1]}\times \mathrm{Hom}(L',M)_{Y}\rightarrow
\bigcup_{[L]}\mathrm{Hom}(L',M\oplus X)^{Y,Z[1]}_{L}
$$
\end{Prop}
\begin{proof}
The additivity of the Hom functor shows there is an isomorphism
$$
\mathrm{Hom}(X,L)\times \mathrm{Hom}(M,L)\cong
\mathrm{Hom}(M\oplus X,L)
$$
for any $X,M,L\in \mc.$ This induces a map $$
\mathrm{Hom}(X,L)_{Y}\times \mathrm{Hom}(M,L)_{Z[1]}\rightarrow
\bigcup_{[L']}\mathrm{Hom}(M\oplus X,L)^{Y,Z[1]}_{L'[1]}
$$
It is a bijection simply since
$$
\mathrm{Hom}(M\oplus
X,L)^{Y,Z[1]}_{L'[1]}=(\mathrm{Hom}(X,L)_Y\times
\mathrm{Hom}(M,L)_{Z[1]})\cap \mathrm{Hom}(M\oplus X,L)_{L'[1]}
$$
It is similar to prove the second result.
\end{proof}

\begin{Thm} The $\bbq$-space $\mh$
is an associative algebra with the $\bbq$-basis $\{u_{[X]}\mid
X\in\mc\},$ the multiplication: $$
u_{[X]}*u_{[Y]}=\sum_{[L]}F_{XY}^{L}u_{[L]}
$$ where $F_{XY}^{L}=\{X,Y\}\cdot\sum_{\alpha\in
V(X,Y;L)}\frac{|\mathrm{End} X_1(\alpha)|}{|\mathrm{Aut}
X_1(\alpha)|}$ and the unit $u_0.$
\end{Thm}
\begin{proof}
For $X,Y,Z \mbox{ and }M\in\mc,$ we need to prove
$u_Z*(u_X*u_Y)=(u_Z*u_X)*u_Y.$ It is equivalent to prove
$$
\sum_{[L]}F_{XY}^{L}F_{ZL}^{M}=\sum_{[L']}F_{ZX}^{L'}F_{L'Y}^{M}.
$$We know that
$$
\sum_{[L]}F_{XY}^{L}F_{ZL}^{M}=\sum_{[L]}\frac{|\mathrm{Hom}(X,L)_Y|}{|\mathrm{Aut}(X)|\{X,X\}}\{X,L\}\cdot
\frac{|\mathrm{Hom}(M,L)_{Z[1]}|}{|\mathrm{Aut}(L)|\{L,L\}}\{M,L\}
$$
By Proposition~\ref{9}, it equals
$$
\frac{1}{|\aut
X|\cdot\{X,X\}}\sum_{[L]}\sum_{[L']}\frac{|\mathrm{Hom}(M\oplus
X,L)^{Y,Z[1]}_{L'[1]}|}{|\aut L|}\cdot \frac{\{M\oplus
X,L\}}{\{L,L\}}
$$
Dually, the right hand side  is
$$
\frac{1}{|\aut
X|\cdot\{X,X\}}\sum_{[L']}\sum_{[L]}\frac{|\mathrm{Hom}(L',M\oplus
X)^{Y,Z[1]}_{L}|}{|\aut L'|}\cdot \frac{\{L',M\oplus
X\}}{\{L',L'\}}
$$
By Proposition~\ref{5}, they equal to each other. We finish the
proof of the theorem.
\end{proof}

\end{document}